\IfFileExists{asl.cls}{\documentclass[bsl]{asl}}{\documentclass[11pt]{amsart}}
\usepackage[T1]{fontenc}
\usepackage[utf8]{inputenc}
\usepackage{lmodern,microtype,amsmath,amssymb,mathtools,enumitem,booktabs}
\usepackage[hidelinks]{hyperref}
\newcommand{\func}[1]{\operatorname{#1}}
\IfFileExists{asl.cls}{}{\usepackage[margin=0.9
in]{geometry}}

\newcommand{\Th}{\operatorname{Th}}

\newtheorem{theorem}{Theorem}[section]
\newtheorem{proposition}[theorem]{Proposition}
\newtheorem{corollary}[theorem]{Corollary}

\IfFileExists{asl.cls}{\address{Department of Philosophy, Linguistics and Theory of Science, University of Gothenburg, Gothenburg, Sweden}}{}

\makeatletter
\IfFileExists{asl.cls}{}{%
  \renewcommand{\@setauthors}{%
    \begingroup
    \centering
    \vspace{5pt}
    \noindent
    \begin{minipage}[t]{0.88\textwidth}
      \centering
      \footnotesize\scshape Ali Enayat\par
      
      \normalfont\scriptsize
      Department of Philosophy, Linguistics, and Theory of Science\par
      University of Gothenburg\par
      ali.enayat@gu.se\par
    \end{minipage}\par
    \vspace{10pt}
    \endgroup
  }%
}
\makeatother

\makeatletter
\IfFileExists{asl.cls}{}{%
\renewcommand{\section}{\@startsection{section}{1}%
  \z@{.7\linespacing\@plus\linespacing}{.5\linespacing}%
  {\normalfont\bfseries\centering}}%
\renewcommand{\subsection}{\@startsection{subsection}{2}%
  \z@{.6\linespacing\@plus.7\linespacing}{.35\linespacing}%
  {\normalfont\scshape\centering}}%
}
\makeatother

\begin{document}

\title[Lessan's Models of Arithmetic]{Hamid Lessan's \\
\emph{Models of Arithmetic}}

\author{Ali Enayat}

\date{\today}

\begin{abstract}
We revisit Hamid Lessan's 1978 Manchester dissertation, \emph{Models of Arithmetic}, an important contribution to the model theory of arithmetic. Surveying its principal results on saturation, definability, initial segments, and self-embeddings, we highlight the mathematical ideas underlying them and discuss their reception and connections with subsequent developments, bringing into focus both the lasting impact of Lessan's work and his lesser-known contributions. In this way, the retrospective also offers a
broader perspective on several central themes
in the model theory of Peano arithmetic and
their development over the past half-century.
\end{abstract}

\maketitle

\tableofcontents

\clearpage

\section{INTRODUCTION}

Hamid Lessan's 1978 doctoral dissertation, \emph{Models of Arithmetic}, completed at the University of Manchester, contains a substantial collection of results concerning nonstandard models of arithmetic \cite{Lessan1978}.%
\footnote{%
The spelling \emph{Lesan} also appears in the literature.}
Lessan's academic career came to an abrupt end under dramatic circumstances, detailed in Asghari's colorful portrait \cite{Asghari2026}. Consequently, he did not have the opportunity to publish his dissertation results or develop them further. Thanks to the efforts of Costas Dimitracopoulos and his co-editors of the volume \cite{WeakArithmetics2013}, the dissertation was republished in 2013.

It is instructive to compare Lessan's dissertation with Alex Wilkie's \emph{Models of Number Theory} (1973) \cite{Wilkie1973}, written at Bedford College, University of London, under the supervision of Wilfrid Hodges. The two dissertations share several central concerns, notably initial segments, self-embeddings, definability, and extensions of models of arithmetic, but differ in emphasis. Wilkie's work focuses particularly on elementary substructures, initial self-embeddings, and end extension problems, whereas Lessan gives greater prominence to saturation and coding methods, including their applications to fragments of arithmetic. The dissertations thus offer complementary perspectives on problems central to the model theory of arithmetic in the 1970s.

Although several of Lessan's findings became known to specialists and influenced subsequent research, other contributions remained relatively obscure. Our purpose in preparing this retrospective is to examine the mathematical significance and historical context of his principal results, trace their connections with subsequent developments, and draw attention to contributions that have received little notice, even among specialists. In doing so, we find that Lessan's dissertation provides a natural vantage point from which to survey some of the central themes of the model theory of Peano arithmetic and the interplay of ideas that has driven their development over the past half-century.

Lessan's dissertation belongs to the tradition of research in mathematical logic that developed at the University of Manchester during the 1970s, particularly in the model theory of arithmetic and the study of fragments of Peano arithmetic. Under the leadership of Jeff Paris and George Wilmers, Manchester became a major centre for research in these areas, attracting a remarkable group of doctoral students---including Laurie Kirby, Hamid Lessan, Costas Dimitracopoulos, Alan Woods, and Richard Kaye---whose work helped shape the subsequent development of the subject. The aforementioned volume \cite{WeakArithmetics2013} also includes Alan Woods's previously unpublished 1981 dissertation.

A recurring theme of Lessan's dissertation is the interaction among saturation, definability, and the structure of nonstandard models of arithmetic. Its principal contributions include methods for constructing saturated end extensions with prescribed theories using arithmetized completeness, the use of partial satisfaction predicates and overspill to obtain recursively saturated structures, and results concerning definable substructures, automorphisms, and reducts. Among the latter are a partial answer to a problem of Gaifman regarding arithmetically presented models with nonarithmetical complete theories, as well as a characterization of countable nonstandard models of Presburger arithmetic that admit expansions to models of Peano arithmetic. Another significant achievement is his characterization of proper initial self-embeddability of models of arithmetic.

In what follows, we present a selection of Lessan's principal results, followed by a discussion of their mathematical significance and connections with subsequent developments. We include proof sketches selectively, concentrating on results that are central to Lessan's work or have had a lasting influence, as well as lesser-known results whose proofs contain particularly illuminating ideas. Our aim is to convey the underlying methods without attempting a systematic reconstruction of the proofs in the dissertation. While every effort has been made to maintain
objectivity, the author's sensibilities and interests have
undoubtedly shaped the presentation.

\begin{quote}
\noindent\textbf{Convention on numbering.}
We use our own numbering for the results stated in this paper, while indicating Lessan's original numbering in parentheses in their headings. Other references to Lessan's results use the numbering in his dissertation.
\end{quote}

Lessan's original bibliography, reproduced in
the Appendix, offers a glimpse of the literature
that informed his research, including unpublished
dissertations, lecture notes, and manuscripts
that circulated among researchers in the model
theory of arithmetic during the 1970s.

\section*{Author's Note}

I learned about Lessan's thesis when I was a
graduate student in the early 1980s, thanks to
spellbinding expositions by C.~Smory\'nski,
especially \cite{Smorynski1981}. Over the
subsequent decades, as I became more knowledgeable
about nonstandard models of arithmetic and had
the opportunity to read and understand his thesis,
the idea of writing a survey of this type on
Lessan's work took root in my mind. The decisive
impetus arrived in the form of Amir Asghari's
biographical account of Lessan \cite{Asghari2026},
and the subsequent companion article
\cite{BahramiEnayat2026}, co-authored with
Saeideh Bahrami, who provided much-needed
camaraderie and feedback. I am also grateful to
Roman Kossak, Costas Dimitracopoulos, Lawrence Wong, and
Leszek Ko\l odziejczyk for their constructive comments
on earlier drafts.

I relied on technical assistance from ChatGPT for typesetting, locating certain references, and editorial matters. The mathematical content, scholarly judgments, and final decisions concerning the presentation are my own.

\section{PRELIMINARIES}

We first fix the arithmetic theories and logical complexity classes used below, then describe cuts, extensions, and embeddings, and finally recall standard systems and saturation. $\mathbb{N}$ is used both for the set of natural numbers, and for the standard model of
arithmetic, the intended meaning being clear from context. Similarly, $M$ is used to denote both a model
of arithmetic and its underlying set when no confusion is likely.

\subsection{Fragments of arithmetic}

We work in the first-order language of arithmetic $\mathcal{L}%
_{A}=\{0,S,+,\cdot ,<\}$. Let $\mathrm{Q}$ be Robinson arithmetic, the familiar finitely
axiomatized theory (see Tarski, Mostowski, and Robinson
\cite{TarskiMostowskiRobinson1953}). We take the ordering in
Robinson arithmetic $\mathrm{Q}$ to be defined by $x<y\iff \exists z\,(z\neq
0\wedge x+z=y)$; thus our use of $<$ is a definitional extension of the
usual language of $\mathrm{Q}$. Peano arithmetic ($\mathrm{PA}$) is the
theory obtained by adjoining to $\mathrm{Q}$ the full induction scheme 
\[
\bigl(\varphi (0,\bar{y})\wedge \forall x\,(\varphi (x,\bar{y})\rightarrow
\varphi (Sx,\bar{y}))\bigr)\rightarrow \forall x\,\varphi (x,\bar{y}),
\]%
for every $\mathcal{L}_{A}$-formula $\varphi (x,\bar{y})$.

For each $n\in \mathbb{N}$, let $\overline{n}$ denote the numeral for $n$.
Every model $M$ of $\mathrm{Q}$ contains the canonical initial segment $\{%
\overline{n}^{M}:n\in \mathbb{N}\}$ of the interpretations of the numerals. This initial segment is the longest well-founded initial segment
of $M$, and isomorphic to $\mathbb{N}$. We identify this canonical copy with 
$\mathbb{N}$, regarding $\mathbb{N}$ as an initial segment of $M$. For our
purposes, an $\mathcal{L}_{A}$-structure $M$ satisfying  $\mathrm{Q}$ is
called \emph{standard} if $M$ is isomorphic to $\mathbb{N}$, and \emph{%
nonstandard} otherwise.

We use the usual arithmetical hierarchy: $\Delta _{0}=\Sigma _{0}=\Pi _{0}$
consists of formulae with bounded quantifiers, and $\Sigma _{n+1}$
(respectively $\Pi _{n+1}$) consists, up to the usual equivalences, of
existential (respectively universal) blocks followed by $\Pi _{n}$
(respectively $\Sigma _{n}$) formulae. We write $I\Gamma $ for induction
restricted to formulae in $\Gamma $, and $B\Sigma _{n}$ for the $\Sigma _{n}$
collection scheme.

The theory $\mathrm{PA}_{\Pi_n}$ consists of
the $\Pi_n$-consequences of $\mathrm{PA}$.
Lessan denotes this theory by $B_n$ and
writes $B$ for $B_2$. We consistently use
$\mathrm{PA}_{\Pi_n}$ in place of this
notation.

Given a model $M$ of arithmetic and a class $\Gamma$ of formulae, we use $\func{Th}_{\Gamma}(M)$ to denote:
\[
\{\sigma : \sigma \text{ is a $\Gamma$-sentence and } M\models\sigma\}.
\]
In particular, $\func{Th}(M)$ denotes the complete theory of $M$.

We use $\mathrm{Exp}$ for the assertion that exponentiation is total; note that $\mathrm{Exp}$ is provable in $%
I\Sigma _{1}$, but not provable in $I\Delta _{0}$. 

\subsection{Cuts, extensions, and embeddings}

An \emph{initial segment} of a model $M$ is a downward-closed subset of its
underlying order. A \emph{cut} of $M$ is a nonempty initial segment of $M$
with no greatest element; when closed under the arithmetic operations it
carries the induced structure. We say that $N$ \emph{end extends}
$M$ if $M$ is a submodel of $N$ whose underlying set is an initial
segment of $N$, and write $M\subseteq_{\mathrm e}N$. If the
extension is also elementary, we write $M\preccurlyeq_{\mathrm e}N$.
An extension $N$ of $M$ is \emph{cofinal} if every element of $N$
is below some element of $M$. More generally, $M\preccurlyeq _{\Gamma }N$
denotes elementarity with respect to formulae in $\Gamma $. An embedding $%
j:M\rightarrow N$ is \emph{initial} when $j[M]$ is an initial segment of $N$%
; an \emph{initial self-embedding} has $N=M$. We call it \emph{proper} when
its image is a proper initial segment.

For models of $\mathrm{PA}$, we follow Lessan's prime-divisibility
convention (Definition~1.3.12 of the dissertation). Write $p_x$ for
the $x$-th prime, with the indexing interpreted internally in a model
$M$. For $a\in M$ and an initial segment $I$ of $M$, put
\[
X_a^I=\{x\in I:M\models p_x\mid a\}.
\]
We say that $X\subseteq I$ is \emph{coded in $M$} if $X=X_a^I$ for
some $a\in M$. Lessan explicitly gives the definition for
$I=\mathbb N$, calling the resulting family of subsets of
$\mathbb N$ the \emph{reals of $M$}, denoted $R_M$; the extension
to arbitrary initial segments is our convention. For a cut $I$ of $M$, its \emph{relative standard system} is
\[
\func{SSy}_I(M)=\{X_a^I:a\in M\}.
\]
In particular, the usual \emph{standard system} is
\[
\func{SSy}(M)=\func{SSy}_{\mathbb N}(M)=R_M.
\]

 Standard systems serve as invariants when comparing models and
their initial segments. An important property of standard systems is their
invariance under end extensions of nonstandard models: if
$M\subseteq_{\mathrm e}N$ are models of
$\mathrm{PA}$ (or even $I\Delta_0$) and $M$ is nonstandard, then:
\[
\func{SSy}(M)=\func{SSy}(N).
\]

\subsection{Definability and saturation}

A \emph{type} $\Phi(\bar x)$ over a finite tuple $\bar a$ in a model $M$
is a finitely satisfiable set of formulae with free variables $\bar x$
and parameters from $\bar a$. The type is \emph{recursive} if its
formulae admit a recursive (computable) enumeration. The model $M$ is
\emph{recursively saturated} if every finitely satisfiable recursive
type over finitely many parameters is realized in $M$.

The notion of recursive saturation, now a standard concept of model theory,
was itself a recent innovation when Lessan undertook his dissertation. It
emerged in the mid-1970s through the work of Barwise and Schlipf \cite%
{BarwiseSchlipf1975,BarwiseSchlipf1976}, alongside independent work of
Ressayre \cite{Ressayre1977} on closely related saturation and resplendence phenomena. Wilmers \cite{Wilmers1980} developed the related notion of
$\mathcal{A}$-saturation, which we will encounter in the next paragraph. Lessan's dissertation belongs to this
formative period. For contemporaneous developments concerning recursively
saturated and cofinal extensions, see also \cite%
{Smorynski1978Abstract,SmorynskiStavi1980,Kotlarski1983}.
Smory\'nski's influential expositions, especially his account of
recursively saturated models \cite{Smorynski1981} and his broader
historical lectures on nonstandard models of arithmetic
\cite{Smorynski1984}, helped organize and disseminate the theory
as it developed.

The dissertation
employs Wilmers's notion of $\mathcal{A}$-saturation, where $\mathcal{A}%
\subseteq \mathcal{P}(\mathbb{N})$ is a \emph{Scott set}: it is closed under
Turing reducibility and join, and contains a path through each infinite
binary tree coded in $\mathcal{A}$. Roughly speaking, $\mathcal{A}$%
-saturation requires realization of consistent partial types coded in $%
\mathcal{A}$, together with appropriate coding of realized complete types
(see Definition~1.4.2 of the dissertation). For $\mathcal{A}$-saturated
models of $\mathrm{PA}$, Lessan proves that $\mathcal{A}=\func{SSy}(M)$
(Lessan's Lemma~2.1.1 \cite{Lessan1978}). Thus the standard system records the information available
for the saturation constructions. Moreover, for a model $M$ that satisfies a modicum of $\mathrm{PA}$, recursive saturation is equivalent to
$\func{SSy}(M)$-saturation. Thus, in the setting considered here,  recursive
saturation and $\func{SSy}(M)$-saturation can be used interchangeably.

For a recursively saturated model $M\models\mathrm{PA}$, the complete
types over finite tuples realized in $M$ are precisely those whose sets
of G\"odel numbers belong to $\func{SSy}(M)$. This is the type-coding
characterization associated with Wilmers's $\mathcal A$-saturation
\cite{Wilmers1980}, together with Lessan's identification
$\mathcal A=\func{SSy}(M)$ (Lessan's Lemma~2.1.1 \cite{Lessan1978}).

\section{END EXTENSIONS AND INITIAL SEGMENTS}

\label{sec:saturation}

We begin with Lessan's constructions of end extensions satisfying prescribed theories and levels of saturation, and compare them with related end-extension criteria. We then turn to his constructions of initial segments within a fixed model and their connections with later work.

\subsection{End extensions}

Lessan's Chapter~2 develops constructions of $\mathcal A$-saturated
end extensions. A principal tool in these constructions is the
arithmetized completeness theorem, in the form recorded by Lessan as
Theorem~2.2.1 of the dissertation, following McAloon. It produces end
extensions equipped with satisfaction relations whose behavior on coded
formulae can be expressed in the original model. In the proof of
Theorem~2.2.2, Lessan uses this internal satisfaction relation to express
the finite satisfiability of coded types; overspill then yields their
realization in the extension. Thus arithmetized completeness provides not
only the extension itself but also the mechanism underlying its
saturation. A representative result is the following.

\begin{theorem}[Lessan, Theorem~2.2.2
  \protect\cite{Lessan1978}]
Let $M$ be a nonstandard model of $\mathrm{PA}$, let
$n\in\mathbb N$, and let $T\supseteq\mathrm{PA}$ be a theory
coded in $\func{SSy}(M)$. If all $\Pi_{n+1}$ consequences
of $T$ hold in $M$, then $M$ has a $\Sigma_n$-elementary
end extension that is a recursively saturated model
of $T$.
\end{theorem}

A further comparison result of Lessan concerns two starting models rather
than a prescribed target theory.

\begin{proposition}[Lessan, Proposition~2.2.7
  \protect\cite{Lessan1978}]
Let $M_1,M_2$ be countable nonstandard models of $\mathrm{PA}$,
and let $n\in\mathbb N$. $M_1$ and $M_2$ admit isomorphic $\Sigma_n$-elementary
end extensions if and only if
\[
\operatorname{SSy}(M_1)=\operatorname{SSy}(M_2)
\]
and their $\Pi_{n+1}$-theories are comparable under inclusion.
\end{proposition}

Lessan also gives an end-extension criterion for \emph{complete
theories} (Proposition~2.2.20 of the dissertation \cite{Lessan1978}). Its coding hypothesis
requires that $\func{Th}_{\Pi_m}(T)$ belong to
$\func{SSy}(M)$ for every $m\in\mathbb N$, together
with the corresponding $\Pi_{n+1}$ compatibility
condition. His proof uses a sequence of end extensions satisfying
successively stronger fragments of the prescribed complete theory.
At about the same time, Wilkie \cite{Wilkie1977}  obtained a related
characterization for complete theories. We state it in
the hierarchical form given by Schmerl
\cite{Schmerl1992}. 

\begin{theorem}[Wilkie--Schmerl
  \protect\cite{Wilkie1977,Schmerl1992}] \label{Wilkie-Schmerl}
Let $M$ be a nonstandard model of $\mathrm{PA}$,
let $n\in\mathbb N$, and let $T\supseteq\mathrm{PA}$
be a complete theory. Then $M$ has a
$\Sigma_n$-elementary end extension satisfying $T$
if and only if
\[
\func{Rep}(T)\subseteq\func{SSy}(M)
\quad\text{and}\quad
\func{Th}_{\Pi_{n+1}}(T)
\subseteq\func{Th}_{\Pi_{n+1}}(M),
\]
where $\func{Rep}(T)$ is the standard system of
the prime model of $T$.
\end{theorem}

Lessan ends Chapter~2 by asking for necessary and sufficient
conditions for a nonstandard model of $\mathrm{PA}$ to have a recursively saturated elementary end extension. This question is closely connected with the
subsequent study of \emph{short recursive saturation}.
Kossak \cite{Kossak1983} characterized short recursively
saturated models in terms of their occurrence as elementary
initial segments of recursively saturated models, thereby
answering Lessan's question. Shochat \cite{Shochat2006}
subsequently developed the theory further, particularly
its connections with expansions and automorphism groups.

\subsection{Initial segments}

Besides constructing end extensions, Lessan studied
initial segments of a fixed model that satisfy prescribed
theories or retain specified degrees of elementarity.
The following theorem illustrates how these requirements
can be imposed simultaneously.

\begin{theorem}[Lessan, Theorem~2.2.9
  \protect\cite{Lessan1978}]\label{thm:lessan-229}
Let $n\in\mathbb N$, and suppose that $M$ is a
countable nonstandard model of $\mathrm{PA}$ with
$\mathbb N\preccurlyeq_{\Sigma_n}M$.
Let $\Phi\subseteq\Th(\mathbb N)$ be coded in $M$.
Then, for every nonstandard $a\in M$, there is a
nonstandard proper initial segment $M'$ of $M$
such that
\[
M'\preccurlyeq_{\Sigma_n}M,\qquad
a\notin M',\qquad M'\models \Phi.
\]
\end{theorem}

Thus the initial segment can be chosen to omit any
prescribed nonstandard element while satisfying the
coded theory $\Phi$. Theorem~\ref{thm:lessan-229} also plays a role
in Lessan's proof of Theorem~\ref{thm:lessan-definable-coinitial}, discussed in
Section~5. A related result of Lessan is the following.
\begin{proposition}[Lessan, Corollary~5.2.2(ii)
  \protect\cite{Lessan1978}]\label{prop:lessan-522}
Every countable nonstandard model $M$ of
$\mathrm{PA}$ has arbitrarily large initial
segments that are themselves models of
$\mathrm{PA}$ and are recursively saturated.
\end{proposition}
Other parts of Lessan's Corollary~5.2.2 \cite{Lessan1978}
strengthen this conclusion. In particular, if
$\operatorname{Th}(M)\in\operatorname{SSy}(M)$, the recursively saturated
initial segments can be chosen elementarily equivalent to $M$;
other clauses give conditions for obtaining arbitrarily large
initial segments satisfying a prescribed complete extension of $\mathrm{PA}$.

The existence of initial segments satisfying
$\mathrm{PA}$ was also investigated by McAloon
\cite{McAloon1982}, who proved that every
countable nonstandard model of $I\Delta_0$
has a nonstandard initial segment satisfying
$\mathrm{PA}$. D'Aquino \cite{DAquino1993}
extended this investigation to arbitrary
models of $I\Delta_0$. Her main theorem characterizes, in terms of the Wainer hierarchy below $\varepsilon_0$, those elements $a$ that belong to an initial segment satisfying $\mathrm{PA}$. For each fixed $a$, the values $f_\alpha(a)$, as $\alpha$ ranges over the ordinals below $\varepsilon_0$, must exist in the ambient model and have a common upper bound there.
D'Aquino's proof uses recursive saturation
of bounded initial structures, obtained
through an internal truth definition.
This provides another connection with
Lessan's bounded-satisfaction construction,
discussed in Section~4. D'Aquino notes that
the result can also be obtained by the venerable method
of indicators. This method was developed by Paris
\cite{Paris1980}, and is presented by Kaye
\cite{Kaye1991}.

Lessan also investigated the approximation of
initial submodels by initial segments satisfying
stronger arithmetic theories. A useful tool is his
\emph{Immortality of the Submodels} principle, first
established for models of $\mathrm{PA}$.

\begin{proposition}[Lessan, Proposition~2.2.19
  (Immortality of the Submodels) \protect\cite{Lessan1978}]
Let $M_1\preccurlyeq_{\Sigma_n}M_2$ be countable
models of $\mathrm{PA}$, where $n\in\mathbb N$.
Let $a,b\in M_1$ and let $T$ be a theory coded in
$\operatorname{SSy}(M_1)$. Then the following are equivalent:
\begin{enumerate}[label=\textup{(\roman*)}]
\item There is a $\Sigma_n$-elementary initial
submodel $I_1$ of $M_1$ satisfying $T$ such that
$a\in I_1$ and $b\notin I_1$.
\item There is a $\Sigma_n$-elementary initial
submodel $I_2$ of $M_2$ satisfying $T$ such that
$a\in I_2$ and $b\notin I_2$.
\end{enumerate}
\end{proposition}

The proof expresses the existence of such an initial
submodel by a coded infinitary sentence. Lessan then
uses saturated end extensions and an isomorphism
argument to transfer the assertion between the two
ambient models. The coding assumption on $T$ makes
the infinitary requirements available in the relevant
standard system. Later in the dissertation, Lessan presented a related transfer
principle for models of the weaker theory
$\mathrm{PA}_{\Pi_2}$. Here the ambient models need
only be related by inclusion, but the conclusion no
longer requires the initial submodels to be
$\Sigma_n$-elementary.

\begin{proposition}[Lessan, Proposition~4.1.18
  (Immortality of the Submodels) \protect\cite{Lessan1978}]
Let $M_1\subseteq M_2$ be countable models of
$\mathrm{PA}_{\Pi_2}$, let $a,b\in M_1$, and let
$T$ be a theory coded in $M_1$. Then the following
are equivalent:
\begin{enumerate}[label=\textup{(\roman*)}]
\item There is an initial submodel $I_1$ of $M_1$
satisfying $T$ such that $a\in I_1$ and $b\notin I_1$.
\item There is an initial submodel $I_2$ of $M_2$
satisfying $T$ such that $a\in I_2$ and $b\notin I_2$.
\end{enumerate}
\end{proposition}

Lessan derives this version by adapting the earlier
argument, using his Lemma~4.1.17 to compare bounded
relational structures in the weaker arithmetic setting.
It is this second form of Immortality that enters the proof of the following theorem, which Lessan notes was independently established by Paris.

\begin{theorem}[Lessan, Theorem~5.2.3
  \protect\cite{Lessan1978}]
If $M_1$ is an initial submodel of $M_2$ and both
models satisfy $\mathrm{PA}_{\Pi_2}$, then
$M_1$ is either the union or the intersection
of a sequence of initial segments of $M_2$
that are models of $\mathrm{PA}$.
\end{theorem}

\begin{proof}[Outline of Lessan's proof]
We describe the argument when $M_2$ is countable (Lessan leaves the uncounable case to the reader).
Fix $b\in M_1$ and $a\in M_2\setminus M_1$.
By Lessan's Lemma~3.1.3, there is a countable
model $M_3\models\mathrm{PA}$ such that
$M_1\preccurlyeq_{\Sigma_1}M_3$. A compactness argument produces a countable elementary
extension $M_4$ of $M_2$ containing $M_3$ and satisfying
$a>c$ for every $c\in M_3$. Indeed, the existential closure of any finite part
of the atomic diagram of $M_3$, with parameters from
$M_1$, is a $\Sigma_1$-formula. Since
$M_1\preccurlyeq_{\Sigma_1}M_3$, the corresponding
finite diagram can be realized in $M_1$, and hence
below $a$ in $M_2$.
Let $\overline M_3$ be the initial segment of $M_4$
generated downward by $M_3$. Since $M_3$ is cofinal in
$\overline M_3$, Gaifman's theorem gives
$M_3\preccurlyeq\overline M_3$, so
$\overline M_3\models\mathrm{PA}$. Moreover,
$b\in\overline M_3$ and $a\notin\overline M_3$.
Applying Immortality to $M_2\preccurlyeq M_4$
with $T=\mathrm{PA}$, we obtain an initial segment
$I$ of $M_2$ satisfying $\mathrm{PA}$ such that
$b\in I$ and $a\notin I$.

Thus $\mathrm{PA}$-initial segments separate any
element of $M_1$ from any element of
$M_2\setminus M_1$. If the union of the $\mathrm{PA}$-initial
segments contained in $M_1$ is cofinal in $M_1$,
that union equals $M_1$. Otherwise, the separation
property implies that the $\mathrm{PA}$-initial
segments containing $M_1$ have intersection $M_1$.
Countability allows the approximating families to
be chosen as sequences. 
\end{proof}

In Remarks~5.2.4, Lessan gives examples showing that each of the two alternatives in the theorem can hold without the other.

\section{LOCAL SATISFACTION}
\label{sec:bounded-satisfaction}

Lessan's Chapter~4 develops a satisfaction
predicate for bounded relational structures.
The construction is noteworthy because it
provides a uniform definition of satisfaction
for all standard formulae, using sufficiently
large codes for finite truth-value computations.
It also illustrates how coding methods can
yield saturation properties in structures
that are coded in nonstandard models of
arithmetic.

\subsection{Local satisfaction}

For a model $K$ of arithmetic and
$\alpha\in K$, let $I_\alpha^K$ denote
the relational structure induced on
$[0,\alpha)^K$, with addition and
multiplication represented by their graphs.
$I_\alpha^K$ is sometimes denoted as $K\upharpoonright\alpha$;
the notation $[0,x]^K$ corresponds to
$I_{x+1}^K$.

Suppose that $K$ is a nonstandard model of $\mathrm{PA}$, $\alpha\in K$ is nonstandard,
and
\[
J\preccurlyeq I_\alpha^K.
\]
The elementarity hypothesis allows coding constructions
carried out in $I_\alpha^K$ to be
transferred to $J$. Lessan uses iterated exponentiation,
defined internally in $K$ by
\[
E_0(x)=x,\qquad
E_{i+1}(x)=2^{E_i(x)}.
\]
An element $\beta\in J$ is called
\emph{small relative to the ambient interval} if
\[
E_\gamma(\beta)<\alpha
\]
for some nonstandard $\gamma\in K$.
In particular, $E_m(\beta)<\alpha$
for every standard $m$. The significance
of this condition is that it provides
sufficient coding space for the
truth-value computations associated
with formulae of arbitrary standard
complexity.

\begin{theorem}[Lessan, Lemma~4.1.4
  \protect\cite{Lessan1978}]
There is a formula $S(u,v,w)$ in the
relational language of arithmetic such
that, whenever $\beta\in J$ is small,
for every standard relational formula
$\varphi(\bar x)$ and every tuple
$\bar a\in(I_\beta^J)^{<\omega}$,
\[
I_\beta^J\models\varphi(\bar a)
\quad\Longleftrightarrow\quad
J\models S\bigl(
\beta,\ulcorner\varphi\urcorner,
\langle\bar a\rangle\bigr).
\]
\end{theorem}

\begin{proof}[Sketch of proof]
The objective is to construct a single formula expressing the
existence of a coded evaluation of the formula whose G\"odel
number is supplied as a parameter. Fix a standard formula $\varphi$.
Let $q$ be the number of its subformulae
and let $r$ bound the number of variables
needed in their evaluation. A truth table
for a subformula with at most $r$ free
variables has at most $\beta^r$ entries.
The truth tables for all subformulae
can therefore be concatenated into a
binary string of length at most
$q(\beta+1)^r$, and coded by a number
below
\[
2^{q(\beta+1)^r}.
\]

The correctness of such a code is
expressed uniformly by the Tarski
recursion clauses. Atomic formulae are first evaluated in the ambient interval
$I_\beta^K$, using its restricted arithmetic relations; the
corresponding evaluations in $I_\beta^J$ are obtained by
elementarity and relativization, as explained below. Boolean connectives
are handled by the corresponding
operations on truth values; and
existential quantification is evaluated
by taking the disjunction over all
assignments below $\beta$.
These conditions can be expressed by
a single first-order formula in the
relational language, using bounded
coding of sequences and computations.

For each fixed standard $\varphi$,
there is a standard integer
$m=m(\varphi)$ such that the truth
tables and the auxiliary computations
needed to verify their correctness
can be coded below $E_m(\beta)$.
The smallness hypothesis ensures
that this bound lies below $\alpha$.
Consequently, the ambient structure
$I_\alpha^K$ contains a code for a
correct evaluation of $\varphi$. Since $J\preccurlyeq I_\alpha^K$,
the existence of such a code transfers
to $J$. Moreover, for $\beta\in J$,
relativization of quantifiers to
$x<\beta$ shows that
\[
I_\beta^J\preccurlyeq I_\beta^K.
\]
Indeed, the interpretation of a relational formula in the interval
below $\beta$ is given by relativizing its quantifiers to that
interval in the ambient structure. Elementarity of $J$ in
$I_\alpha^K$ therefore gives the asserted elementarity.
An external induction on the standard
subformulae of $\varphi$ verifies that
the transferred evaluation assigns
the correct truth values in
$I_\beta^J$. The assertion that the
evaluation code assigns value $1$
to $\varphi(\bar a)$ defines the
required formula $S$.
\end{proof}
The generality of allowing an elementary substructure $J\preccurlyeq I_\alpha^K$ is used in Lessan's subsequent arguments, notably in transferring the recursive-saturation conclusions to bounded intervals of models of $\mathrm{PA}_{\Pi_1}$. We return to this application in Section 4.4.

\subsection{Bounded satisfaction}

The construction above has a useful
formulation over the weak theory
$I\Delta_0$. We follow the formulation
from \cite[Theorem~2.2]{EnayatWong2017}. For $s\in M\models I\Delta_0$, we write
\[
2^{s^{\mathbb N}}\subseteq_{\mathrm e}M
\]
if the elements below some $2^{s^k}$,
for standard $k$, form a proper initial
segment of $M$. In particular, all
these powers exist in $M$ and have a
common upper bound.

Similarly, the notation
$t>2^{s^{\mathbb N}}$ means
\[
t>2^{s^k}
\qquad\text{for every standard }
k\in\mathbb N.
\]

\begin{theorem}[Bounded satisfaction;
Lessan \protect\cite{Lessan1978}]
There is a $\Delta_0$-formula
$\operatorname{Sat}_0(u,v,w)$ such that,
whenever $M\models I\Delta_0$ and
$s,t\in M$ satisfy
\[
t>2^{s^{\mathbb N}},
\]
then, for every standard formula
$\varphi(\bar x)$ of the relational
language of arithmetic and every
$\bar a\in(I_s^M)^{<\omega}$,
\[
I_s^M\models\varphi(\bar a)
\quad\Longleftrightarrow\quad
M\models\operatorname{Sat}_0
(\ulcorner\varphi\urcorner,
 \langle\bar a\rangle,t).
\]
\end{theorem}

The bound on $t$ is uniform over
all standard formulae. It supplies
room for coding the truth-value
computations required by formulae
of arbitrarily large standard
complexity. 

Lessan's work on local satisfaction is a key ingredient in the proof of the theorem below, which was first established by Paris and Dimitracopoulos \cite{ParisDimitracopoulos1982}.

\begin{theorem}[Truth for bounded formulae]
There is a $\Sigma_1$-formula
$\func{Tr}_{\Delta_0}(u,v)$ such that,
for every $\Delta_0$-formula
$\varphi(\bar x)$,
$I\Delta_0+\mathrm{Exp}$ proves
\[
\forall\bar x\,
\bigl(
\func{Tr}_{\Delta_0}
(\ulcorner\varphi\urcorner,
 [\bar x])
\leftrightarrow
\varphi(\bar x)
\bigr).
\]
\end{theorem}

The $\Sigma_1$-complexity arises
by existentially quantifying over
an auxiliary bound for a coded
truth-value computation.
H\'ajek and Pudl\'ak  use this
construction to establish the
finite axiomatizability of
$I\Delta_0+\mathrm{Exp}$
\cite[Theorem~5.6]{HajekPudlak1993}.

Adamowicz, Ko\l odziejczyk, and Paris \cite{AdamowiczEtAl2012} subsequently investigated truth definitions in models where exponentiation is not total. They formulated an optimality principle, called \emph{LIO} (Lessan's bound Is Optimal), asserting that the size of the auxiliary parameter in Lessan's $\Delta_0$ universal formula cannot, in a precise sense, be uniformly reduced. They proved that LIO implies the consistency of
\[
I\Delta_0+\neg\mathrm{Exp}+\neg B\Sigma_1.
\]
The argument exploits the fact that LIO prevents certain sets of true bounded formulae from being coded, whereas $\Sigma_1$-collection would provide the uniform bounds needed for such coding. Thus an optimality question arising from Lessan's truth-definition construction is connected with the failure of collection in weak arithmetic.

\subsection{Recursive saturation}

The local satisfaction theorem also yields
recursive saturation of suitable
initial structures.

\begin{proposition}[Lessan,
Proposition~4.1.5
\protect\cite{Lessan1978}]
Let $K\models\mathrm{PA}$ be nonstandard, let $\alpha\in K$
be nonstandard, and suppose $J\preccurlyeq I_\alpha^K$.
If $\beta\in J$ is small in the sense of Section~4.1,
then $I_\beta^J$ is recursively saturated.
\end{proposition}

The passage from bounded satisfaction
to recursive saturation uses overspill
applied to the satisfiability of coded
finite fragments of recursive types.
The smallness hypothesis supplies the
uniform satisfaction predicate needed
for this argument, which we spell out
in the following corollary.

The following corollary is implicit in Lessan's bounded-satisfaction construction and its recursive-saturation consequence (Lemma~4.1.4 and Proposition~4.1.5 of the dissertation).

\begin{corollary}[
Lessan \protect\cite{Lessan1978}]
Let $M\models I\Delta_0$ and
$s\in M$. If
\[
2^{s^{\mathbb N}}
\subseteq_{\mathrm e}M,
\]

then the induced relational
structure $I_s^M$ is recursively
saturated.
\end{corollary}

\begin{proof}[Sketch of proof]
Let $p(x,\bar a)$ be a recursive type
finitely satisfiable in $I_s^M$, and
write
\[
p=\{\varphi_i(x,\bar a):i\in\mathbb N\}
\]
for a recursive enumeration. We may
assume that $M$ is nonstandard.
Choose a nonstandard $t\in M$ such
that $t>2^{s^k}$ for every standard $k$. By Theorem~4.2, the predicate
$\operatorname{Sat}_0$, with parameter
$t$, correctly expresses satisfaction
in $I_s^M$ for every standard formula.
Fix a Turing machine $e$ computing the
G\"odel numbers of the formulae $\varphi_i$,
and let $C_e(i,c,u)$ be a $\Delta_0$-predicate
expressing that $c$ codes a terminating
computation of $e$ on input $i$ with output $u$.
Consider the $\Delta_0$-formula
\[
\begin{aligned}
\theta(k,t,\bar a)\;:\Longleftrightarrow\;
\exists x<s\,\forall i<k\,\exists c,u<t\,
\bigl(&C_e(i,c,u)\\
&{}\wedge\operatorname{Sat}_0
(u,\langle x,\bar a\rangle,t)\bigr).
\end{aligned}
\]
For each standard $k$, the computations
involved have standard codes, hence codes
below the nonstandard element $t$.
Finite satisfiability of $p$, together
with the correctness of $\operatorname{Sat}_0$
for standard formulae, therefore gives
$\theta(k,t,\bar a)$. By $\Delta_0$-overspill,
$\theta(k,t,\bar a)$ holds for some
nonstandard $k$. A corresponding witness
$b<s$ satisfies every standard formula
of $p$, since the computation predicate
gives the correct G\"odel number for each
standard input and $\operatorname{Sat}_0$
correctly evaluates the resulting formula.
Thus $b$ realizes $p$.
\end{proof}

A later application of the same
underlying method appears in the
model theory of nonstandard finite
structures. Ajtai
\cite{Ajtai2007}
established a completeness theorem
for the existence of suitable
expansions and end extensions
of such structures. In a
subsequent proof, Garl\'{\i}k
\cite{Garlik2015}
uses a bounded truth definition
and overspill to obtain a
relative saturation property
for an internally coded
structure. This is essentially
the mechanism underlying
Lessan's bounded-satisfaction lemma and
its recursive-saturation consequence, although
the surrounding extension
problem is quite different.

\subsection{An application to end extensions}

We conclude with an application of the preceding results
that illustrates the usefulness of allowing an arbitrary
elementary substructure $J$ in Lessan's local satisfaction
theorem. For this application, Lessan uses a related notion
of smallness, applying to an entire initial segment rather
than to an individual element. Following Definition~4.1.3(ii)
of the dissertation \cite{Lessan1978}, an initial segment
$M$ of a model $M'\models\mathrm{PA}_{\Pi_1}$ is said to be
\emph{small in $M'$} if there are $b\in M'\setminus M$
and a nonstandard $\gamma\in M'$ such that $2^{b^\gamma}$
exists in $M'$.

\begin{theorem}[Lessan, Theorem~4.1.11
  \protect\cite{Lessan1978}]
Let $M$ be a countable structure in the language
of arithmetic. Suppose that $M$ has an end
extension $M'\models\mathrm{PA}_{\Pi_1}$
such that $M$ is small in $M'$.
Then $M$ has an end extension satisfying
$\mathrm{PA}$.
\end{theorem}

\begin{proof}[Outline of the proof]
Suppose that $M'$ is countable. Choose $\beta\in M'$ such that
$M\subseteq I_\beta^{M'}$ and
$I_\beta^{M'}$ is small in $M'$.
Lessan's Lemma~4.1.9 shows that
bounded relational intervals of
models of $\mathrm{PA}_{\Pi_1}$
can be elementarily embedded into
corresponding intervals of models
of $\mathrm{PA}$. Thus there is
$K\models\mathrm{PA}$ such that
\[
I_\beta^{M'}\preccurlyeq I_\beta^K.
\]
The recursive-saturation consequence
of the local satisfaction theorem,
discussed above, now shows that
$I_\beta^{M'}$ is recursively saturated. Since this structure is countable,
the Barwise--Schlipf theorem implies
that it is resplendent. In Lessan's
terminology, it is a
$\Sigma(\mathrm{rec})$-structure:
recursive requirements involving
additional relation symbols can be
realized by expansions whenever
they are consistent with the
elementary diagram of the structure. To obtain the required end extension,
Lessan formulates recursive infinitary
requirements expressing the existence
of an extension satisfying
$\mathrm{PA}$ in which $M$ is an
initial segment. Their consistency
is witnessed by the corresponding
bounded interval of $K$.
Resplendence then allows these
requirements to be realized over
$I_\beta^{M'}$, yielding an end
extension of $M$ satisfying
$\mathrm{PA}$. This is the expansion
argument developed in Proposition~4.1.7
of the dissertation.

Now suppose that $M'$ is uncountable.
By the smallness hypothesis, there
are $b\in M'\setminus M$ and a
nonstandard $\gamma\in M'$ such
that $2^{b^\gamma}$ exists in $M'$.
Since $M$ is countable, the downward
L\"owenheim--Skolem theorem gives a
countable elementary submodel
$N\preccurlyeq M'$ containing $M$,
$b$, and $\gamma$. Elementarity ensures that
$N\models\mathrm{PA}_{\Pi_1}$ and
that $2^{b^\gamma}$ exists in $N$.
Moreover, since
$M\subseteq_{\mathrm e}M'$ and
$M\subseteq N\subseteq M'$, we have
$M\subseteq_{\mathrm e}N$.
Thus $M$ remains small in $N$,
and the countable-case argument
applies to $N$, giving the desired
end extension of $M$.
\end{proof}

\section{DEFINABILITY AND REDUCTS}

\label{sec:definability}

This section examines Lessan's results on substructures generated by definable elements and the definability of the standard cut. We then discuss his partial answer to Gaifman's problem on arithmetical presentations and his criteria for expanding reducts to models of Peano arithmetic.

\subsection{Definable hulls and the standard cut}

Let $M\models\mathrm{PA}$. Write $K^n(M;\bar a)$ for the substructure of $M$ consisting of elements $\Sigma_n$-definable over $\bar a$, and $I^n(M;\bar a)$ for the initial segment of $M$ generated by $K^n(M;\bar a)$. We omit $\bar a$ when no parameters are allowed. Lessan developed these constructions systematically in Chapter~3, \S3.1 of his dissertation, establishing their elementary-substructure properties and results concerning induction and collection. 

Paris and Kirby \cite{ParisKirby1978} independently established the corresponding separation phenomena: their Proposition~7 uses the parameter-free definable hull $K^n(M)$ (which they denote by $K_n$) to separate $I\Sigma_{n-1}$ from $B\Sigma_n$, while Proposition~8 uses an increasing-chain construction to separate $B\Sigma_n$ from $I\Sigma_n$. They explicitly acknowledge Lessan's independent proof of the latter result. The following theorem summarizes Lessan's results.
\begin{theorem}[Definable hulls and their convex hulls]\label{definable-hulls}
Let $M\models\mathrm{PA}$, $n\geq1$, and let $\bar a$ be a finite
tuple from $M$ (possibly empty). 
\begin{enumerate}[label=\textup{(\alph*)}]
\item $K^n(M;\bar a)\preccurlyeq_{\Sigma_n}M$.
If $K^n(M;\bar a)$ is nonstandard, then:
\[
 K^n(M;\bar a)\models I\Sigma_{n-1}+\neg B\Sigma_n.
\]
\item $I^n(M;\bar a)\preccurlyeq_{\Sigma_{n-1}}M$.
If $I^n(M;\bar a)$ is nonstandard, then:
\[
 I^n(M;\bar a)\models B\Sigma_n+\neg I\Sigma_n.
\]
\end{enumerate}
\end{theorem}

Lessan's Proposition~3.1.6 gives further structural information
about the two parameter-free hull constructions.

\begin{proposition}[Lessan, Proposition~3.1.6
  \protect\cite{Lessan1978}]\label{prop:lessan-316}
Let $M,N\models\mathrm{PA}$ and $n\geq1$.
\begin{enumerate}[label=\textup{(\alph*)}]
\item $K^n(M)\cong K^n(N)$ if and only if
$\operatorname{Th}_{\Sigma_n}(M)=\operatorname{Th}_{\Sigma_n}(N)$.
\item Every element of $K^n(M)$ is definable without parameters
in $K^n(M)$; in particular, this structure is rigid.
\item If $M,N$ are countable and $I^n(M)$ is nonstandard, then:
$I^n(M)\cong I^n(N)$ if and only if
\[
\operatorname{Th}_{\Sigma_n}(M)=\operatorname{Th}_{\Sigma_n}(N)
\quad\text{and}\quad
\operatorname{SSy}(M)=\operatorname{SSy}(N).
\]
\item If $M$ is countable and $I^n(M)$ is nonstandard, then
$I^n(M)$ has $2^{\aleph_0}$ automorphisms.
\end{enumerate}
\end{proposition}

The contrast between the two hull constructions is particularly
striking: passing from the parameter-free definable hull to the
initial segment it generates can radically change the automorphism
structure. Lessan's isomorphism criterion for initial-segment hulls
uses his comparison theorem for saturated end extensions
(Proposition~2.2.7 of the dissertation). The abundance of
automorphisms was subsequently established for the broader class
of countable models of $B\Sigma_n+\mathrm{Exp}+\neg I\Sigma_n$
by Kossak \cite{Kossak1990} and Kaye \cite{Kaye1991}.
For parameterized hulls, every automorphism of $K^n(M;\bar a)$
fixing $\bar a$ pointwise is the identity.

There is also a separation between successive definability levels.
Lessan's Theorem~3.1.15 \cite{Lessan1978} shows that, for $n\geq1$,
if $K^{n+1}(M)$ is nonstandard, then $K^n(M)$ is not cofinal
in $K^{n+1}(M)$. Lessan also investigates a complementary
phenomenon: whether increasing definitional complexity can yield
nonstandard elements below all the nonstandard elements definable
at the preceding level. Write $A\downarrow B$ when some
nonstandard element of $B$ lies below every nonstandard element of $A$.

\begin{theorem}[Lessan, Theorem~4.2.8
  \protect\cite{Lessan1978}]\label{thm:lessan-428}
For every $n\in\mathbb N$ there is a model $M\models\mathrm{PA}$
such that
\[
K^n(M)\downarrow K^{n+1}(M).
\]
\end{theorem}

The proof uses coding at successive levels of the arithmetical
hierarchy, drawing on an independent-family construction of
Mostowski and Friedman. The downward-separation phenomenon is
also relevant to Lessan's subsequent investigation of initial
segments satisfying true low-complexity arithmetic; see
Corollary~5.2.10 of the dissertation.

\begin{theorem}[Lessan, Corollary~3.1.11
  \protect\cite{Lessan1978}]\label{thm:lessan-standard-cut}
Let $M\models\mathrm{PA}$ and let $n\geq1$.
Then:
\begin{enumerate}[label=\textup{(\alph*)}]
\item The standard cut $\mathbb N$ is definable without parameters
in $K^n(M)$ and $I^n(M)$.
\item The complete theories of $K^n(M)$ and $I^n(M)$ are
nonarithmetical.
\item $\operatorname{Th}(K^n(M))$ is Turing reducible to
$\operatorname{Th}(K^{n+1}(M))$.
\end{enumerate}
\end{theorem}

The second conclusion follows by relativizing arithmetic sentences
to the definable standard cut; the third expresses a further
computational relationship between successive definable hulls.

Lessan also obtains restrictions on the elementary extensions of
parameter-free definable hulls.

\begin{theorem}[Lessan, Theorem~3.1.17
  \protect\cite{Lessan1978}]\label{thm:lessan-3117}
Suppose that $K^n(M)$ is nonstandard, where $M\models\mathrm{PA}$
and $n\geq1$. If $K^n(M)\prec N$ , then:
\begin{enumerate}[label=\textup{(\alph*)}]
\item $N$ contains a nonstandard initial segment elementarily
equivalent to $\mathbb N$.
\item Between any two elements of $K^n(M)$ that are not
finitely apart, there is an element of $N\setminus K^n(M)$.
\end{enumerate}
Consequently, $K^n(M)$ has neither a proper elementary cofinal
extension nor a proper elementary end extension.
\end{theorem}
Around the same time, McAloon \cite{McAloon1978} used arithmetized completeness and incompleteness arguments to construct models of fragments of Peano arithmetic in which the standard cut is definable at prescribed levels of the arithmetical hierarchy. In the same paper, he also showed that every
recursively enumerable extension of
$\mathrm{PA}$ has a nonstandard model with
a downward cofinal sequence of nonstandard
elements that are both existentially and
universally definable.

A related result on definability of the standard cut was obtained
by Adamowicz, Cord\'on-Franco, and Lara-Mart\'in
\cite{AdamowiczEtAl2016}: the standard cut is parameter-free
$\Pi_1$-definable in every existentially closed model of
$I\Delta_0+\mathrm{Exp}$. They also show that this complexity bound
is optimal. Unlike Lessan's result, this concerns existentially
closed models rather than definable hulls; the comparison illustrates
the recurring role of definability of the standard cut in the study
of fragments of arithmetic.

A more recent application of Lessan's definable-hull results concerns interpretability.
Enayat and \L e\l yk \cite{EnayatLelyk2024}
use Theorems~\ref{definable-hulls} and~\ref{thm:lessan-standard-cut} to construct distinct bi-interpretable
extensions of fragments of arithmetic. This establishes failures of
\emph{tightness}, in contrast with $\mathrm{PA}$, which is tight by a
theorem of Visser; see Enayat and \L e\l yk
\cite{EnayatLelyk2024} for further details. Gruza, Ko\l odziejczyk, and \L e\l yk \cite{GruzaKolodziejczykLelyk2025} have also used the definable-hull construction in their investigation of tightness and solidity in fragments of arithmetic. Their results show, in particular, that for every $n$ there are proper subtheories of $\mathrm{PA}$ extending $I\Sigma_n$ that are solid, as well as ones that are tight but not neat; the solid theories can moreover be chosen not to interpret $\mathrm{PA}$.

We return now to Lessan's investigation of definable elements,
turning from their structural properties to their distribution
within the ambient model. 

\begin{theorem}[Lessan, Theorem~2.2.17 \protect\cite{Lessan1978}]\label{thm:lessan-definable-coinitial}
Let $n\in\mathbb{N}$, and suppose $M$ is a nonstandard model of $\mathrm{PA}$
with $\mathbb{N}\preccurlyeq_{\Sigma_n}M$. Then the nonstandard
$\Delta_{n+1}$-definable elements of $M$ are coinitial in the
nonstandard part of $M$ if and only if
$\func{Th}_{\Pi_{n+1}}(\mathbb{N})\notin\func{SSy}(M)$.\footnote{Note that in the published 2013 reprint of Lessan's dissertation,
the statement of Theorem~2.2.17 has
$\func{Th}_{\Pi_{n+1}}(\mathbb N)\in\func{SSy}(M)$.
The proof of the theorem, as well as Corollary~2.2.18, shows that
the correct condition is
$\func{Th}_{\Pi_{n+1}}(\mathbb N)\notin\func{SSy}(M)$.}
\end{theorem}

Lessan notes that this theorem answers a problem of McAloon
\cite{McAloon1978}, who had independently constructed models
exhibiting the coinitiality phenomenon. He also credits Jensen
and Ehrenfeucht \cite{JensenEhrenfeucht1976} with an independent
solution for $n=0$. The Jensen-Ehrenfeucht paper also developed, independently and in a
different setting, ideas concerning recursive saturation and the
construction of models of prescribed theories. As noted by Smory\'nski \cite{Smorynski1984}, several results that subsequently became
part of the theory of recursively saturated models were already present
in Jensen and Ehrenfeucht's work.

The distribution of definable elements also plays a role in controlling the ranges of initial self-embeddings. Mijajlovi\'c \cite{Mijajlovic1983}, attributing the result to Marker and Wilkie, showed that if $M$ is a countable model of $\mathrm{PA}$ and no parameter-free $\Delta_1$-definable element of $M$ lies above $a\in M$, then there is an initial self-embedding of $M$ whose image omits $a$.\footnote{Here, the hypothesis is equivalent to requiring that every true parameter-free $\Sigma_1$ sentence of $M$ have a witness below $a$. The result may therefore be viewed as a ``target version'' of Friedman's self-embedding theorem; compare Kaye \cite[Theorems~12.3 and~12.4]{Kaye1991}. We are grateful to Leszek Ko\l odziejczyk for pointing out this equivalent formulation.} 

\subsection{A problem of Gaifman}

A model of arithmetic is called \emph{arithmetical} if it has a presentation
on a subset of $\mathbb{N}$ with arithmetical operations and relations. The standard model of PA is the canonical example of such a model, and by Tarski's undefinability of truth, its theory is not arithmetical. This prompted Gaifman \cite{Gaifman1972} to ask: \emph{Is there a nonstandard arithmetical model of $\mathrm{%
PA}$ whose complete theory is nonarithmetical?} Lessan obtains the following partial answer. His construction starts with an arithmetically presented model of $\mathrm{PA}+\neg\mathrm{Con}(\mathrm{PA})$ and passes to a suitable parameter-free definable hull. Lemma~3.1.4(i) \cite{Lessan1978} supplies the required restricted elementarity, Proposition~3.2.2 \cite{Lessan1978} preserves arithmetical presentability, and the definability of the standard cut (Theorem~\ref{thm:lessan-standard-cut}), together with Lessan's preceding analysis of definability and complete theories, yields the required nonarithmeticality.

\begin{theorem}[Lessan, Theorem~3.2.4 \protect\cite{Lessan1978}]
For every $n\geq1$ there is a nonstandard arithmetical model $M\models%
\mathrm{PA}_{\Pi_n}$ such that $\func{Th}(M)$ is not arithmetical.
\end{theorem}

Lessan also observes (following Proposition~3.2.2 of the
dissertation) that if $M\models\mathrm{PA}$ is arithmetically
presented, then, for every $n\geq1$,
\[
\operatorname{Th}(K^n(M))\equiv_T\operatorname{Th}(\mathbb N)
\equiv_T\operatorname{Th}(I^n(M)).
\]
Thus the nonarithmeticality of the complete theories of these
hulls can be sharpened to an exact Turing-degree assertion.

Gaifman's problem was eventually answered
affirmatively by Leo Harrington, who constructed a nonstandard arithmetical
model of $\mathrm{PA}$ whose complete theory is not arithmetical. The construction used the method of infinitely many workers.
Knight \cite{Knight1986} gives a related account and strengthened the
result by obtaining such a model of low Turing degree.

\subsection{Reducts and expansions}

A \emph{reduct} forgets the interpretations of symbols from a structure's language, while an \emph{%
expansion} interprets additional symbols on the same underlying set. A
different question addressed in the dissertation is when a reduct can be
expanded to a model of $\mathrm{PA}$. Lessan establishes a general criterion
for this problem.

In the result below, we follow Lessan's terminology. Let $L'$ be a sublanguage of $L$, let $T$ be an $L$-theory, and let $T'$ be its reduct to $L'$. The theory $T'$ is said to be \emph{complete in $T$} if there is an $L$-formula $\tau(u,v)$ such that, for every $L'$-formula $\varphi(\bar x)$,
\[
T\vdash\forall\bar x\,\bigl(\tau(\ulcorner\varphi\urcorner,\langle\bar x\rangle)\leftrightarrow\varphi(\bar x)\bigr).
\]
The theory $T'$ is \emph{strongly complete in $T$} if it is complete in the usual syntactic sense and is complete in $T$ in Lessan's sense \cite[Definitions~2.3.1 and~2.3.3]{Lessan1978}.

Lessan also proves a uniqueness result for reducts
\cite[Corollary~2.3.4]{Lessan1978}. In the arithmetic setting
of his Chapter~2, suppose that $T'$ is strongly complete in
$T$ and that $M_1,M_2\models T$ have the same standard system.
Then their reducts to $L'$ are $L_{\infty\omega}$-equivalent;
if $M_1$ and $M_2$ are countable, the reducts are isomorphic.
The result does not assume that the given models are saturated.
Thus Lessan's analysis of reducts yields both uniqueness results
and criteria for expandability.

\begin{theorem}[Lessan, Corollary~2.3.6 \protect\cite{Lessan1978}]
Let $L^{\prime }$ be a sublanguage of the language of arithmetic, and let $%
\mathrm{PA}^{\prime }$ be a reduct of $\mathrm{PA}$ that is strongly
complete in $\mathrm{PA}$. A countable nonstandard model
of $\mathrm{PA}^{\prime }$ admits an expansion to a model of $\mathrm{PA}$
if and only if it is recursively saturated.
\end{theorem}

A particularly transparent application concerns the additive language. Recall that \emph{%
Presburger arithmetic} is the complete first-order theory of the additive
structure of $\mathbb{N}$.

The following characterization was obtained by Lessan \cite{Lessan1978} and
independently by Jensen and Ehrenfeucht \cite{JensenEhrenfeucht1976}, as reported by Wilmers \cite{Wilmers1985}.

\begin{corollary}[Expandability of Presburger models]
A countable nonstandard model of Presburger arithmetic admits an expansion
to a model of $\mathrm{PA}$ if and only if it is recursively saturated.
\end{corollary}

Subsequent work developed the result in two directions. C\'egielski, McAloon, and Wilmers \cite{CegielskiMcAloonWilmers1982} showed that full Peano arithmetic can be replaced by the much weaker theory $I\Delta_0$: a countable nonstandard model of Presburger arithmetic is recursively saturated if and only if it admits an expansion to a model of bounded induction. 
The countability assumption of the corollary, on the other hand, cannot
simply be dropped. Lipshitz and Nadel \cite{LipshitzNadel1978} proved that recursive saturation
remains necessary for an additive reduct to admit a PA-expansion, but is
insufficient for uncountable models. Knight \cite{Knight1983} subsequently
considered expansion to a fixed completion $T$ of $\mathrm{PA}$ in
the uncountable setting. For Presburger models of cardinality
$\aleph_1$, she characterized when an $L_{\infty\omega}$-equivalent
model admits an expansion to a model of $T$, using recursive
saturation together with coding conditions on the $\Pi_n$
consequences of $T$. The conclusion concerns an
$L_{\infty\omega}$-equivalent model, not necessarily an expansion
of the given model itself. Harnik  \cite%
{Harnik1986} subsequently showed that a countable recursively saturated
Presburger model has many nonisomorphic $\omega _{1}$-like recursively
saturated elementary end extensions, while only one of their isomorphism
types admits an expansion to a model of PA. 

\section{SELF-EMBEDDINGS AND AUTOMORPHISMS}
\label{sec:comparison}

The study of self-embeddings and automorphisms of models of arithmetic
has undergone extensive development since Lessan's dissertation.
Considerably more subsequent work is closely related to his results
in these areas than to most of the other topics considered here; this
accounts for the greater space devoted to later developments. We treat
the two subjects separately, placing the relevant refinements and
analogues alongside Lessan's results.

\subsection{Initial self-embeddings and comparison theorems}

The historical point of departure is Friedman
\cite{Friedman1973}, who proved that every countable nonstandard
model of $\mathrm{PA}$ is isomorphic to a proper initial segment
of itself. Wilkie \cite{Wilkie1977} subsequently obtained the
following comparison theorem.

\begin{theorem}[Wilkie]\label{Wilkie-comparison}
Let $M$ and $N$ be countable nonstandard models of $\mathrm{PA}$.
Then $M$ is isomorphic to arbitrarily high proper initial segments
of $N$ if and only if
\[
\func{SSy}(M)=\func{SSy}(N)\quad\text{and}\quad
\func{Th}_{\Pi_2}(M)\subseteq\func{Th}_{\Pi_2}(N).
\]
\end{theorem}

The result stated below as Theorem~\ref{thm:lessan-514} gives a related
initial-segment comparison criterion under
weaker arithmetic hypotheses. Following Lessan \cite[Definition~5.1.2]{Lessan1978}, a model $M$ is
\emph{1-tall} if its parameter-free $\Sigma_1$-definable elements are not
cofinal in $M$. It is \emph{1-extendable} if it has an end extension
$M^*$ satisfying $I\Delta_0$ such that
$\func{Th}_{\Sigma_1}(M)=\func{Th}_{\Sigma_1}(M^*)$.
These conditions concern, respectively, definable witnesses and the
existence of a suitable end extension.

\begin{theorem}[Lessan, Proposition~5.1.4 \protect\cite{Lessan1978}]\label{thm:lessan-514}
Suppose $M_1$ and $M_2$ are countable structures in the language of
arithmetic, $M_1$ is 1-extendable, and $M_2\models\mathrm{PA}_{\Pi_2}$ is
1-tall. Then $M_1$ is isomorphic to an initial segment of $M_2$ if and only if 
\[
\func{SSy}(M_1)=\func{SSy}(M_2)  \quad\text{and}\quad  \func{Th}%
_{\Sigma_1}(M_1)\subseteq\func{Th}_{\Sigma_1}(M_2). 
\]
\end{theorem}

Theorem~\ref{thm:lessan-514} retains standard systems and low-complexity
truths as invariants while weakening the arithmetic assumptions. Specializing to a single model gives the following corollary.

\begin{theorem}[Lessan, Corollary~5.1.5 \protect\cite{Lessan1978}]
A countable model $M\models\mathrm{PA}_{\Pi_2}$ is isomorphic to a proper
initial segment of itself if and only if it is 1-tall and 1-extendable.
\end{theorem}

Lessan further identifies a canonical initial segment from the
ranges of initial self-embeddings.

\begin{corollary}[Lessan, Corollary~5.1.8
  \protect\cite{Lessan1978}]\label{cor:lessan-518}
If a countable model $M\models\mathrm{PA}_{\Pi_2}$ is isomorphic
to an initial segment of itself, then:
\[
I^1(M)=\bigcap\{I\subseteq_{\mathrm e}M:I\cong M\}.
\]
\end{corollary}

Thus the initial segment generated by the parameter-free
$\Sigma_1$-definable elements can also be recovered as the
intersection of the ranges of initial self-embeddings. This
concerns the images of embeddings, rather than their fixed points.

Subsequent work refined the theory of initial
self-embeddings in two directions: weakening
the arithmetic hypotheses and obtaining
greater control over the embeddings.
H\'ajek and Pudl\'ak
\cite{HajekPudlak1980} established a
parameter-preserving initial-embedding theorem
under hypotheses involving a shared cut,
coded subsets, and $\Sigma_1$-theories.
Dimitracopoulos \cite{Dimitracopoulos1985}
obtained a version of Friedman's theorem
over $I\Delta_0+B\Sigma_2$.
Dimitracopoulos and Paris
\cite{DimitracopoulosParis1988}, and
independently Ressayre \cite{Ressayre1987},
extended Friedman's theorem to $I\Sigma_1$.
Ressayre also obtained the following
characterization, which allows prescribed
initial segments to be fixed pointwise.

\begin{theorem}
For a countable nonstandard model
$M\models I\Delta_0$, the following
are equivalent:
\begin{enumerate}[label=\textup{(\alph*)}]
\item $M\models I\Sigma_1$;
\item For every $a\in M$, there is a
proper initial self-embedding of $M$
fixing every $x\leq a$.
\end{enumerate}
\end{theorem}

In a related direction, Dimitracopoulos and Paris \cite{DimitracopoulosParis1988}
extended Lessan's result (Theorem~\ref{thm:lessan-514}) to a setting involving
$I\Delta_0+B\Sigma_1+\exp$. As summarized by Bahrami
\cite{Bahrami2024}, a parameterized formulation for countable
nonstandard models of $I\Delta_0+\mathrm{Exp}$ continues to use standard
systems and $\Sigma_1$-theories as comparison invariants, with suitable
tallness and extendability assumptions. These later formulations use their own conventions for tallness
and extendability. 

The location of initial images between prescribed
bounds was also investigated by Wilkie. As recorded
by Shavrukov \cite{Shavrukov2016}, Wilkie obtained
a criterion for the existence of an initial segment
$I$ of a model $M\models\mathrm{PA}$ satisfying
$a<I<b$ and $I\equiv M$; when $M$ is countable,
$I$ can moreover be chosen isomorphic to $M$.
Kaye \cite{Kaye1991} gives a systematic treatment
of related parameterized initial-embedding criteria. Lessan ends Chapter~5 with an explicit question (p.~416): \emph{what is the
weakest theory $T\supseteq\mathrm{PA}_{\Pi_1}$ such that every model of $T$
is isomorphic to an initial segment of itself?} It appears that this question remains open.

Beyond prescribing the image of an embedding, one can prescribe its
fixed-point structure. Bahrami and Enayat
\cite{BahramiEnayat2018} characterized the cuts that can occur as the
longest initial segment fixed pointwise, or as the full fixed-point set,
of such an embedding. The full fixed-point set need not coincide with the longest
initial segment fixed pointwise. Bahrami \cite{Bahrami2024} subsequently extended the question from cuts
to small submodels. For a countable nonstandard model of $I\Sigma_1$,
her characterization relates the strongness of a cut $I$ to the property
that every $I$-small $\Sigma_1$-elementary submodel can be realized as
the full fixed-point set of a proper initial self-embedding. 

A key generalization of the aforementioned theorems of Dimitracopoulos--Paris and Ressayre was established by Tanaka in the context of models of the fragment  $\mathrm{WKL}_0$ of second-order arithmetic. 

\begin{theorem}[Tanaka \protect\cite{Tanaka1997}]
Every countable nonstandard model $(M,\mathfrak{X})$ of $\mathrm{WKL}_0$ is isomorphic to a
proper initial part of itself, with an embedding fixing any prescribed
bounded initial region.
\end{theorem}

A subsequent proof of Tanaka's theorem was given in \cite{Enayat2013}. Bahrami
\cite{Bahrami2020} subsequently characterized the cuts that can occur
as the longest initial segment fixed pointwise, or as the full fixed-point
set, of a proper initial self-embedding of a countable nonstandard model
of $\mathrm{WKL}_0$.

\subsection{Automorphisms }

Lessan also investigated the automorphism
groups of nonstandard models of arithmetic
and their initial submodels. His results
establish the existence of continuum many
automorphisms under suitable saturation
hypotheses and, more generally, reveal a
connection between the automorphisms of an
initial segment and the number of its
isomorphic copies within an ambient model.

\begin{proposition}[Lessan, Proposition~2.1.10
  \protect\cite{Lessan1978}]
Let $M$ be a countable $\mathcal A$-saturated
model of $\mathrm{PA}$. Then every nonstandard
proper elementary initial submodel $M'$ of $M$
has $2^{\aleph_0}$ automorphisms.
\end{proposition}

Lessan also observes that every countably
infinite $\mathcal A$-saturated structure
has continuum many automorphisms. It is also routine
to deduce from the results of Barwise and Schlipf
\cite{BarwiseSchlipf1976} that every countably infinite
recursively saturated structure has $2^{\aleph_0}$
automorphisms. Lessan's Theorem~5.3.8 \cite{Lessan1978}, using the notions of
\emph{clear} and \emph{strongly clear}
initial segments, yields a more general
dichotomy.

\begin{corollary}[Lessan, Corollary~5.3.9 \protect\cite{Lessan1978}]
If $I$ is a nonstandard initial segment of a countable model $M$, then
either $I$ has $2^{\aleph_0}$ automorphisms or $M$ has $2^{\aleph_0}$
initial segments isomorphic to $I$.
\end{corollary}

The systematic study of automorphisms of countable recursively saturated
models of $\mathrm{PA}$ developed substantially in the years after Lessan's
dissertation. Beginning with work of Smory\'nski and Kotlarski, and
continuing through the investigations of Kaye, Kossak, Kotlarski, and
Schmerl, attention turned from the mere abundance of automorphisms to the
structure of their groups and the submodels and cuts they can fix;
see Smory\'nski \cite{Smorynski1981}, Kaye, Kossak, and Kotlarski
\cite{KayeKossakKotlarski1991}, and Kossak and Schmerl
\cite{KossakSchmerl2006}. 

Two representative results illustrate the greater precision of this theory.
Smory\'nski \cite{Smorynski1981} showed that every cut closed under exponentiation in a countable
recursively saturated model of $\mathrm{PA}$ can occur as the \emph{longest
initial segment fixed pointwise} by an automorphism. Kaye, Kossak, and Kotlarski \cite{KayeKossakKotlarski1991}
characterized those cuts that can occur as \emph{full fixed-point sets} in
terms of strongness and elementarity. An automorphism may fix isolated elements above its longest pointwise
fixed initial segment. The strongness condition for full fixed-point
sets has a connection with second-order arithmetic: by the
Kirby--Paris characterization \cite{ParisKirby1978}, a cut $I$ of
a model $M$ of $\mathrm{PA}$ is strong if and only if
\[
(I,\func{SSy}_I(M))\models\mathrm{ACA}_0.
\]

A further development concerns prescribed fixed-point substructures.
Enayat \cite{EnayatUnified2007} investigated the realization of
fixed-point structures by automorphisms of recursively saturated models,
settling a conjecture of Schmerl (stated as Question~9 in Kossak and
Schmerl \cite{KossakSchmerl2006}).

\section*{Appendix: Lessan's Bibliography}

To document the sources available to Lessan, we reproduce the bibliography of
Lessan's dissertation as printed in the 2013 reprint
\cite[pp.~417--418]{Lessan1978}, preserving the original order,
reference labels, and bibliographical information. In particular,
Lessan dates Wilkie's dissertation to 1971; its correct date is 1973.

\begingroup
\footnotesize
\setlength{\parindent}{0pt}
\setlength{\parskip}{4pt}
\begin{quote}
\noindent\textbf{[B\&S.]} J. Barwise and J. Schlipf, ``An Introduction to Recursively Saturated and Resplendent Models'', J.S.L. Vol. 41 (1976).\par
\noindent\textbf{[C\&K.]} C. C. Chang and H. J. Keisler, ``Model Theory'', North Holland (1973).\par
\noindent\textbf{[F.]} H. Friedman, ``Countable Models of Set Theories'', Cambridge Summer School in Mathematical Logic, 1971, Springer Lecture Notes 37, pp. 539--573.\par
\noindent\textbf{[G.]} H. Gaifman, ``A Note on Models and Submodels of Arithmetic'', Conference in Mathematical Logic, London'70, Springer Lecture Notes 255, pp. 128--144.\par
\noindent\textbf{[H\&W.]} J. Hirschfeld and W. Wheeler, ``Forcing, Arithmetic, Division Rings'', Springer Lecture Notes 454.\par
\noindent\textbf{[J\&E.]} D. Jensen and A. Ehrenfeucht, ``Some Problems in Elementary Arithmetics'', Fundamenta Mathematicae (1976).\par
\noindent\textbf{[K.]} L. A. S. Kirby, Ph.D. Thesis, Manchester (1977).\par
\noindent\textbf{[Ku.]} D. W. Kueker, ``Back-and-Forth Arguments and Infinitary Logics'', Infinitary Logic: In Memoriam Carol Karp, Springer Lecture Notes 492, pp. 17--73 (1975).\par
\noindent\textbf{[Mc.]} K. McAloon, ``Completeness Theorems, Incompleteness Theorems and Models of Arithmetic'', to appear in TAMS.\par
\noindent\textbf{[P\&K.]} J. B. Paris and L. A. S. Kirby, ``$\Sigma_n$-Collection Schemas in Arithmetic'' (mimeographed).\par
\noindent\textbf{[S.]} D. Scott, ``Algebras of Sets Binumerable in Complete Extensions of Arithmetic'', Proc. Sympos. Pure Math. Vol. 5, A.M.S. (1962).\par
\noindent\textbf{[Sh.]} J. R. Shoenfield, ``Mathematical Logic'', Addison-Wesley (1967).\par
\noindent\textbf{[Sm.]} C. Smorynski, ``Consistence and Related Mathematical Properties'', University of Amsterdam Report 75-02.\par
\noindent\textbf{[Wk.]} A. J. Wilkie, Ph.D. Thesis, Bedford College, London, 1971.\par
\noindent\textbf{[W1.]} G. M. Wilmers, Ph.D. Thesis, Oxford, 1975.\par
\noindent\textbf{[W2.]} G. M. Wilmers, Lecture Notes, Wroclaw, Winter, 1976.\par
\end{quote}
\endgroup


\begingroup

\begin{thebibliography}{99}
\bibitem{AdamowiczEtAl2012} {\footnotesize Z.~Adamowicz, L.~A. Ko\l %
odziejczyk, and J.~B. Paris, \emph{Truth definitions without exponentiation
and the $\Sigma_1$ collection scheme}, J. Symbolic Logic \textbf{77} (2012),
no.~2, 649--655. }

\bibitem{AdamowiczEtAl2016} {\footnotesize Z.~Adamowicz,
A.~Cord\'on-Franco, and F.~F. Lara-Mart\'in,
\emph{Existentially closed models in the framework of arithmetic},
J. Symbolic Logic \textbf{81} (2016), no.~2, 774--788. }

\bibitem{Ajtai2007} {\footnotesize M.~Ajtai, \emph{Generalizations of the
compactness theorem and G\"odel's completeness theorem for nonstandard
finite structures}, in \emph{Computer Science Logic}, Lecture Notes in
Computer Science, vol.~4484, Springer, 2007, pp.~13--33. }

\bibitem{Asghari2026} {\footnotesize A.~Asghari, \emph{Hamid Lesan, a
logician larger than life}, People and Mathematics, July 2026 (English and
Persian versions), %
\url{https://mathswell.com/people-and-mathematics/articles/hamid-lesan.html}%
. }

\bibitem{Bahrami2020} {\footnotesize S.~Bahrami, \emph{Tanaka's theorem
revisited}, Arch. Math. Logic \textbf{59} (2020), 865--877. }

\bibitem{Bahrami2024} {\footnotesize S.~Bahrami, \emph{Self-embeddings of
models of arithmetic; fixed points, small submodels, and extendability}, J.
Symbolic Logic \textbf{89} (2024), no.~3, 1044--1066. }

\bibitem{BahramiEnayat2018} {\footnotesize S.~Bahrami and A.~Enayat, \emph{%
Fixed points of self-embeddings of models of arithmetic}, Ann. Pure Appl.
Logic \textbf{169} (2018), no.~6, 487--513. }

\bibitem{BahramiEnayat2026}
{\footnotesize S.~Bahrami and A.~Enayat,
\emph{Beyond the Natural Numbers},
People and Mathematics, Issue~02, October 2026,
\url{https://mathswell.com/people-and-mathematics/articles/beyond-the-natural-numbers.html}.}

\bibitem{BarwiseSchlipf1975} {\footnotesize J.~Barwise and J.~Schlipf, \emph{%
On recursively saturated models of arithmetic}, in \emph{Model Theory and
Algebra: A Memorial Tribute to Abraham Robinson}, Lecture Notes in
Mathematics, vol.~498, Springer, Berlin, 1975, pp.~42--55. }

\bibitem{BarwiseSchlipf1976} {\footnotesize J.~Barwise and J.~Schlipf, \emph{%
An introduction to recursively saturated and resplendent models}, J.
Symbolic Logic \textbf{41} (1976), no.~2, 531--536. }

\bibitem{CegielskiMcAloonWilmers1982} {\footnotesize P.~C\'egielski,
K.~McAloon, and G.~M. Wilmers, \emph{Mod\`eles r\'ecursivement satur\'es de
l'addition et de la multiplication des entiers naturels}, in \emph{Logic
Colloquium '80}, North-Holland, Amsterdam, 1982, pp.~57--68. }

\bibitem{WeakArithmetics2013} {\footnotesize P.~C\'egielski, C.~Cornaros,
and C.~Dimitracopoulos (eds.), \emph{New Studies in Weak Arithmetics}, CSLI
Lecture Notes, vol.~211, CSLI Publications, Stanford, 2013, \url{https://web.stanford.edu/group/cslipublications/cslipublications/site/9781575867236.shtml}.}

\bibitem{DAquino1993} {\footnotesize P.~D'Aquino, \emph{A sharpened version of McAloon's theorem on initial segments of models of $I\Delta_0$}, Ann. Pure Appl. Logic \textbf{61} (1993), no.~1--2, 49--62. }

\bibitem{Dimitracopoulos1985} {\footnotesize C.~Dimitracopoulos, \emph{A
generalization of a theorem of H. Friedman}, Z. Math. Logik Grundlag. Math. 
\textbf{31} (1985), 221--225. }

\bibitem{DimitracopoulosParis1988} {\footnotesize C.~Dimitracopoulos and
J.~B. Paris, \emph{A note on a theorem of H. Friedman}, Z. Math. Logik
Grundlag. Math. \textbf{34} (1988), no.~1, 13--17. }

\bibitem{EnayatUnified2007} {\footnotesize A.~Enayat, \emph{Automorphisms of
models of arithmetic: A unified view}, Ann. Pure Appl. Logic \textbf{145}
(2007), 16--36. }

\bibitem{Enayat2013} {\footnotesize A.~Enayat, \emph{A new proof of Tanaka's
theorem}, in \emph{New Studies in Weak Arithmetics}, P.~C\'egielski,
C.~Cornaros, and C.~Dimitracopoulos (eds.), CSLI Lecture Notes, vol.~211,
CSLI Publications, Stanford, 2013, pp.~93--102, \url{https://gup.ub.gu.se/file/167548}.}

\bibitem{EnayatLelyk2024} {\footnotesize A.~Enayat and M.~\L e\l yk, \emph{%
Categoricity-like properties in the first order realm}, Journal for the
Philosophy of Mathematics \textbf{1} (2024), 63--98. }

\bibitem{EnayatWong2017} {\footnotesize A.~Enayat and T.~L. Wong, \emph{%
Unifying the model theory of first-order and second-order arithmetic via $%
\mathrm{WKL}_0^*$}, Ann. Pure Appl. Logic \textbf{168} (2017), no.~6,
1247--1283. }

\bibitem{Friedman1973} {\footnotesize H.~M. Friedman, \emph{Countable models
of set theories}, in \emph{Cambridge Summer School in Mathematical Logic},
Lecture Notes in Mathematics, vol.~337, Springer, Berlin, 1973,
pp.~539--573, \url{https://u.osu.edu/friedman.8/foundational-adventures/publications/}.}

\bibitem{GruzaKolodziejczykLelyk2025} {\footnotesize P.~Gruza, L.~A.~Ko\l odziejczyk, and M.~\L e\l yk, \emph{Tightness and solidity in fragments of Peano arithmetic}, arXiv:2512.09120, 2025. \url{https://arxiv.org/abs/2512.09120}. }

\bibitem{Gaifman1972} {\footnotesize H.~Gaifman, \emph{A note on models and
submodels of arithmetic}, in \emph{Conference in Mathematical Logic---London
'70}, Lecture Notes in Mathematics, vol.~255, Springer, Berlin, 1972,
pp.~128--144. }

\bibitem{Garlik2015} {\footnotesize M.~Garl\'ik, \emph{A new proof of
Ajtai's completeness theorem for nonstandard finite structures}, Arch. Math.
Logic \textbf{54} (2015), 413--424. }

\bibitem{HajekPudlak1980} {\footnotesize P.~H\'ajek and P.~Pudl\'ak, \emph{%
Two orderings of the class of all countable models of Peano arithmetic}, in 
\emph{Model Theory of Algebra and Arithmetic}, Lecture Notes in Mathematics,
vol.~834, Springer, 1980, pp.~174--185. }

\bibitem{HajekPudlak1993} {\footnotesize P.~H\'ajek and P.~Pudl\'ak,
\emph{Metamathematics of First-Order Arithmetic}, Perspectives in
Mathematical Logic, Springer, Berlin, 1993. }

\bibitem{Harnik1986} {\footnotesize V.~Harnik, \emph{$\omega_1$-like
recursively saturated models of Presburger's arithmetic}, J. Symbolic Logic 
\textbf{51} (1986), no.~2, 421--429. }

\bibitem{JensenEhrenfeucht1976} {\footnotesize D.~Jensen and A.~Ehrenfeucht, \emph{Some problem in elementary arithmetics}, Fund. Math. \textbf{92} (1976), no.~3, 223--245. }

\bibitem{Kaye1991} {\footnotesize R.~Kaye, \emph{Models of Peano Arithmetic}%
, Oxford Logic Guides, Oxford University Press, Oxford, 1991. }

\bibitem{KayeKossakKotlarski1991} {\footnotesize R.~Kaye, R.~Kossak, and
H.~Kotlarski, \emph{Automorphisms of recursively saturated models of
arithmetic}, Ann. Pure Appl. Logic \textbf{55} (1991), no.~1, 67--99. }

\bibitem{Knight1983} {\footnotesize J.~F. Knight, \emph{Additive structure in uncountable models for a fixed completion of $P$}, J. Symbolic Logic \textbf{48} (1983), no.~3, 623--628. }

\bibitem{Knight1986} {\footnotesize J.~F. Knight, \emph{Effective
construction of models}, in \emph{Logic Colloquium '84}, Studies in Logic
and the Foundations of Mathematics, vol.~120, North-Holland, Amsterdam,
1986, pp.~105--119. }

\bibitem{Kossak1983} {\footnotesize R.~Kossak, \emph{A certain class of
models of Peano arithmetic}, J. Symbolic Logic \textbf{48} (1983), no.~2,
311--320. }

\bibitem{Kossak1990} {\footnotesize R.~Kossak, \emph{On extensions of models
of strong fragments of arithmetic}, Proc. Amer. Math. Soc. \textbf{108}
(1990), no.~1, 223--232. }

\bibitem{KossakSchmerl2006} {\footnotesize R.~Kossak and J.~H. Schmerl, 
\emph{The Structure of Models of Peano Arithmetic}, Oxford Logic Guides,
vol.~50, Clarendon Press, Oxford, 2006, \url{https://doi.org/10.1093/acprof:oso/9780198568278.001.0001}.}

\bibitem{Kotlarski1983} {\footnotesize H.~Kotlarski, \emph{On cofinal
extensions of models of arithmetic}, J. Symbolic Logic \textbf{48} (1983),
no.~2, 253--262. }




\bibitem{Lessan1978} {\footnotesize H.~Lessan, \emph{Models of Arithmetic},
Ph.D. thesis, University of Manchester, February 1978; reprinted in \cite%
{WeakArithmetics2013}, pp.~359--418 (including the prefatory note by
G.~Wilmers), \url{https://web.stanford.edu/group/cslipublications/cslipublications/site/9781575867236.shtml}.}

\bibitem{LipshitzNadel1978} {\footnotesize L.~Lipshitz and M.~Nadel, \emph{%
The additive structure of models of arithmetic}, Proc. Amer. Math. Soc. 
\textbf{68} (1978), no.~3, 331--336. }

\bibitem{McAloon1978} {\footnotesize K.~McAloon, \emph{Diagonal methods and
strong cuts in models of arithmetic}, in \emph{Logic Colloquium '77},
Studies in Logic and the Foundations of Mathematics, vol.~96, North-Holland,
Amsterdam, 1978, pp.~171--181, \url{https://doi.org/10.1016/S0049-237X(08)72000-7}.}

\bibitem{McAloon1978TAMS}
{\footnotesize K.~McAloon,
\emph{Completeness theorems, incompleteness
theorems and models of arithmetic},
Trans. Amer. Math. Soc. \textbf{239}
(1978), 253--277.}

\bibitem{McAloon1982}
{\footnotesize K.~McAloon,
\emph{On the complexity of models of arithmetic},
J. Symbolic Logic \textbf{47} (1982), no.~2,
403--415.}

\bibitem{Mijajlovic1983} {\footnotesize \v{Z}.~Mijajlovi\'c, \emph{Submodels
and definable points in models of Peano arithmetic}, Notre Dame J. Formal
Logic \textbf{24} (1983), no.~4, 417--425. }

\bibitem{Paris1980} {\footnotesize J.~B. Paris, \emph{A hierarchy of cuts in models of arithmetic}, in \emph{Model Theory of Algebra and Arithmetic}, Lecture Notes in Mathematics, vol.~834, Springer, Berlin, 1980, pp.~312--337. }

\bibitem{ParisDimitracopoulos1982} {\footnotesize J.~B. Paris and C.~Dimitracopoulos, \emph{Truth definitions for $\Delta_0$ formulae}, in \emph{Logic and Algorithmic (Zurich, 1980)}, Monographies de l'Enseignement Mathématique, vol.~30, Université de Genève, Geneva, 1982, pp.~317--329. }

\bibitem{ParisKirby1978} {\footnotesize J.~B. Paris and L.~A.~S. Kirby, 
\emph{$\Sigma_n$-collection schemas in arithmetic}, in \emph{Logic
Colloquium '77}, Studies in Logic and the Foundations of Mathematics,
vol.~96, North-Holland, Amsterdam, 1978, pp.~199--209, \url{https://doi.org/10.1016/S0049-237X(08)72003-2}.}

\bibitem{Ressayre1977} {\footnotesize J.-P.~Ressayre, \emph{Models with compactness properties relative to an admissible language}, Ann. Math. Logic \textbf{11} (1977), no.~1, 31--55. }

\bibitem{Ressayre1987} {\footnotesize J.-P. Ressayre, \emph{Nonstandard
universes with strong embeddings, and their finite approximations}, in \emph{%
Logic and Combinatorics}, Contemporary Mathematics, vol.~65, American
Mathematical Society, Providence, RI, 1987, pp.~333--358. }


\bibitem{Schmerl1992} {\footnotesize J.~H. Schmerl, \emph{End extensions of
models of arithmetic}, Notre Dame J. Formal Logic \textbf{33} (1992), no.~2,
216--219. }

\bibitem{Shavrukov2016}
{\footnotesize V.~Yu.~Shavrukov,
\emph{Duality, non-standard elements, and
dynamic properties of r.e.\ sets},
Ann. Pure Appl. Logic \textbf{167}
(2016), 939--981.}

\bibitem{Shochat2006}
{\footnotesize E.~Shochat,
\emph{Countable short recursively saturated models of arithmetic},
Ph.D. thesis, Graduate Center, City University of New York, 2006,
\url{https://academicworks.cuny.edu/gc_etds/1786/}.}

\bibitem{Smorynski1978Abstract} {\footnotesize C.~Smory\'nski, \emph{%
Recursively saturated nonstandard models of arithmetic}, abstract 758-02-2,
Notices Amer. Math. Soc., August 1978. }

\bibitem{Smorynski1981} {\footnotesize C.~Smory\'nski, \emph{Recursively
saturated nonstandard models of arithmetic}, J. Symbolic Logic \textbf{46}
(1981), no.~2, 259--286. }

\bibitem{Smorynski1984} {\footnotesize C.~Smory\'nski,
\emph{Lectures on nonstandard models of arithmetic}, in
\emph{Logic Colloquium '82}, Studies in Logic and the Foundations
of Mathematics, vol.~112, North-Holland, Amsterdam, 1984,
pp.~1--70. }

\bibitem{SmorynskiStavi1980} {\footnotesize C.~Smory\'nski and J.~Stavi, 
\emph{Cofinal extension preserves recursive saturation}, in \emph{Model
Theory of Algebra and Arithmetic}, Lecture Notes in Mathematics, vol.~834,
Springer, Berlin, 1980, pp.~338--345. }

\bibitem{TarskiMostowskiRobinson1953} {\footnotesize A.~Tarski, A.~Mostowski, and R.~M.~Robinson, \emph{Undecidable Theories}, Studies in Logic and the Foundations of Mathematics, North-Holland, Amsterdam, 1953. }

\bibitem{Tanaka1997} {\footnotesize K.~Tanaka, \emph{The self-embedding
theorem of $\mathrm{WKL}_0$ and a non-standard method}, Ann. Pure Appl.
Logic \textbf{84} (1997), no.~1, 41--49. }

\bibitem{Wilkie1973} {\footnotesize A.~J. Wilkie,
\emph{Models of Number Theory}, Ph.D. thesis,
Bedford College, University of London, 1973.}

\bibitem{Wilkie1977} {\footnotesize A.~J. Wilkie, \emph{On the theories of
end-extensions of models of arithmetic}, in \emph{Set Theory and Hierarchy
Theory V}, Lecture Notes in Mathematics, vol.~619, Springer, 1977,
pp.~305--310, \url{https://doi.org/10.1007/BFb0067661}.}

\bibitem{Wilmers1980} {\footnotesize G.~M.~Wilmers, \emph{Minimally saturated models}, in \emph{Model Theory of Algebra and Arithmetic}, Lecture Notes in Mathematics, vol.~834, Springer, Berlin, 1980, pp.~370--380. }

\bibitem{Wilmers1985} {\footnotesize G.~M. Wilmers, \emph{Bounded
existential induction}, J. Symbolic Logic \textbf{50} (1985), no.~1, 72--90. 
}

\end{thebibliography}
\endgroup
\end{document}